\magnification=1200

\hsize=11.25cm    
\vsize=18cm       
\parindent=12pt   \parskip=5pt     

\hoffset=.5cm   
\voffset=.8cm   

\pretolerance=500 \tolerance=1000  \brokenpenalty=5000

\catcode`\@=11

\font\eightrm=cmr8         \font\eighti=cmmi8
\font\eightsy=cmsy8        \font\eightbf=cmbx8
\font\eighttt=cmtt8        \font\eightit=cmti8
\font\eightsl=cmsl8        \font\sixrm=cmr6
\font\sixi=cmmi6           \font\sixsy=cmsy6
\font\sixbf=cmbx6

\font\tengoth=eufm10 
\font\eightgoth=eufm8  
\font\sevengoth=eufm7      
\font\sixgoth=eufm6        \font\fivegoth=eufm5

\skewchar\eighti='177 \skewchar\sixi='177
\skewchar\eightsy='60 \skewchar\sixsy='60

\newfam\gothfam           \newfam\bboardfam

\def\tenpoint{
  \textfont0=\tenrm \scriptfont0=\sevenrm \scriptscriptfont0=\fiverm
  \def\rm{\fam\z@\tenrm}
  \textfont1=\teni  \scriptfont1=\seveni  \scriptscriptfont1=\fivei
  \def\oldstyle{\fam\@ne\teni}\let\old=\oldstyle
  \textfont2=\tensy \scriptfont2=\sevensy \scriptscriptfont2=\fivesy
  \textfont\gothfam=\tengoth \scriptfont\gothfam=\sevengoth
  \scriptscriptfont\gothfam=\fivegoth
  \def\goth{\fam\gothfam\tengoth}
  
  \textfont\itfam=\tenit
  \def\it{\fam\itfam\tenit}
  \textfont\slfam=\tensl
  \def\sl{\fam\slfam\tensl}
  \textfont\bffam=\tenbf \scriptfont\bffam=\sevenbf
  \scriptscriptfont\bffam=\fivebf
  \def\bf{\fam\bffam\tenbf}
  \textfont\ttfam=\tentt
  \def\tt{\fam\ttfam\tentt}
  \abovedisplayskip=12pt plus 3pt minus 9pt
  \belowdisplayskip=\abovedisplayskip
  \abovedisplayshortskip=0pt plus 3pt
  \belowdisplayshortskip=4pt plus 3pt 
  \smallskipamount=3pt plus 1pt minus 1pt
  \medskipamount=6pt plus 2pt minus 2pt
  \bigskipamount=12pt plus 4pt minus 4pt
  \normalbaselineskip=12pt
  \setbox\strutbox=\hbox{\vrule height8.5pt depth3.5pt width0pt}
  \let\bigf@nt=\tenrm       \let\smallf@nt=\sevenrm
  \normalbaselines\rm}

\def\eightpoint{
  \textfont0=\eightrm \scriptfont0=\sixrm \scriptscriptfont0=\fiverm
  \def\rm{\fam\z@\eightrm}
  \textfont1=\eighti  \scriptfont1=\sixi  \scriptscriptfont1=\fivei
  \def\oldstyle{\fam\@ne\eighti}\let\old=\oldstyle
  \textfont2=\eightsy \scriptfont2=\sixsy \scriptscriptfont2=\fivesy
  \textfont\gothfam=\eightgoth \scriptfont\gothfam=\sixgoth
  \scriptscriptfont\gothfam=\fivegoth
  \def\goth{\fam\gothfam\eightgoth}
  
  \textfont\itfam=\eightit
  \def\it{\fam\itfam\eightit}
  \textfont\slfam=\eightsl
  \def\sl{\fam\slfam\eightsl}
  \textfont\bffam=\eightbf \scriptfont\bffam=\sixbf
  \scriptscriptfont\bffam=\fivebf
  \def\bf{\fam\bffam\eightbf}
  \textfont\ttfam=\eighttt
  \def\tt{\fam\ttfam\eighttt}
  \abovedisplayskip=9pt plus 3pt minus 9pt
  \belowdisplayskip=\abovedisplayskip
  \abovedisplayshortskip=0pt plus 3pt
  \belowdisplayshortskip=3pt plus 3pt 
  \smallskipamount=2pt plus 1pt minus 1pt
  \medskipamount=4pt plus 2pt minus 1pt
  \bigskipamount=9pt plus 3pt minus 3pt
  \normalbaselineskip=9pt
  \setbox\strutbox=\hbox{\vrule height7pt depth2pt width0pt}
  \let\bigf@nt=\eightrm     \let\smallf@nt=\sixrm
  \normalbaselines\rm}

\tenpoint

\def\pc#1{\bigf@nt#1\smallf@nt}         \def\pd#1 {{\pc#1} }

\catcode`\;=\active
\def;{\relax\ifhmode\ifdim\lastskip>\z@\unskip\fi
\kern\fontdimen2  -1.2 \fontdimen3 \string;}

\catcode`\:=\active
\def:{\relax\ifhmode\ifdim\lastskip>\z@\unskip\fi\penalty\@M\ \fi\string:}

\catcode`\!=\active
\def!{\relax\ifhmode\ifdim\lastskip>\z@
\unskip\fi\kern\fontdimen2  -1.1 \fontdimen3 \string!}

\catcode`\?=\active
\def?{\relax\ifhmode\ifdim\lastskip>\z@
\unskip\fi\kern\fontdimen2  -1.1 \fontdimen3 \string?}

\frenchspacing

\def\raggedbottom{\topskip 10pt plus 36pt\r@ggedbottomtrue}

\def\pointir{\unskip . --- \ignorespaces}

\def\Medbreak{\vskip-\lastskip\medbreak}

\long\def\th#1 #2\enonce#3\endth{
   \Medbreak\noindent
   {\pc#1} {#2\unskip}\pointir{\it #3}\smallskip}

\def\decale#1{\smallbreak\hskip 28pt\llap{#1}\kern 5pt}
\def\decaledecale#1{\smallbreak\hskip 34pt\llap{#1}\kern 5pt}
\def\puce{\smallbreak\hskip 6pt{$\scriptstyle\bullet$}\kern 5pt}

\def\eqalign#1{\null\,\vcenter{\openup\jot\m@th\ialign{
\strut\hfil$\displaystyle{##}$&$\displaystyle{{}##}$\hfil
&&\quad\strut\hfil$\displaystyle{##}$&$\displaystyle{{}##}$\hfil
\crcr#1\crcr}}\,}

\catcode`\@=12

\showboxbreadth=-1  \showboxdepth=-1

\newcount\numerodesection \numerodesection=1
\def\section#1{\bigbreak
 {\bf\number\numerodesection.\ \ #1}\nobreak\medskip
 \advance\numerodesection by1}

\mathcode`A="7041 \mathcode`B="7042 \mathcode`C="7043 \mathcode`D="7044
\mathcode`E="7045 \mathcode`F="7046 \mathcode`G="7047 \mathcode`H="7048
\mathcode`I="7049 \mathcode`J="704A \mathcode`K="704B \mathcode`L="704C
\mathcode`M="704D \mathcode`N="704E \mathcode`O="704F \mathcode`P="7050
\mathcode`Q="7051 \mathcode`R="7052 \mathcode`S="7053 \mathcode`T="7054
\mathcode`U="7055 \mathcode`V="7056 \mathcode`W="7057 \mathcode`X="7058
\mathcode`Y="7059 \mathcode`Z="705A


\def\hfl#1#2#3{\smash{\mathop{\hbox to#3{\rightarrowfill}}\limits
^{\textstyle#1}_{\textstyle#2}}}

\def\Q{{\bf Q}}

\def\R{{\bf R}}
\def\C{{\bf C}}

\def\Z{{\bf Z}}

\def\F{{\bf F}}

\def\Gal{\mathop{\rm Gal}\nolimits}

\def\to{\rightarrow}

\def\normressym(#1,#2)_#3{\displaystyle\left({#1,#2\over#3}\right)}

\def\mod{\mathop{\rm mod.}\nolimits}

\newcount\refno 
\long\def\ref#1:#2<#3>{                                        
\global\advance\refno by1\par\noindent                              
\llap{[{\bf\number\refno}]\ }{#1} \pointir{\it #2} #3\goodbreak }

\def\citer#1(#2){[{\bf\number#1}\if#2\empty\relax\else,\ {#2}\fi]}

\newbox\bibbox
\setbox\bibbox\vbox{
\bigskip
\centerline{---$*$---$*$---}
\bigbreak
\centerline{{\pc REFERENCE}}

\ref{\pc ARTIN} (E) \& {\pc TATE} (J):
Class field theory,
<AMS Chelsea Publishing, Providence, RI, 2009. viii+194 pp.> 
\newcount\artintate \global\artintate=\refno 

\ref{\pc CASSELS} (J) \& {\pc FR\"OHLICH} (A):
Algebraic number theory,
<Proceedings of the instructional conference held at the University of
Sussex, Brighton, September 1•¡¹17, 1965. Academic Press Inc, London,
1986. xviii+366 pp.>
\newcount\casselsfrohlich \global\casselsfrohlich=\refno

\ref{\pc CHEVALLEY} (C):
Class field theory,
<Nagoya University, Nagoya, 1954. ii+104 pp.>
\newcount\chevalley \global\chevalley=\refno

\ref{\pc DALAWAT} (C):
Six lectures on quadratic reciprocity,
<arXiv:1404.4918.>
\newcount\dalawat \global\dalawat=\refno

\ref{\pc GODEMENT} (R):
Analyse math\'ematique IV, 
<Springer-Verlag, Berlin, 2003. xii+599.>
\newcount\godement \global\godement=\refno

\ref{\pc HASSE} (H):
Bericht \"uber neuere Untersuchungen und Probleme aus der Theorie der
algebraischen Zahlk\"orper. Teil I: Klassenk\"orpertheorie. Teil Ia:
Beweise zu Teil I. Teil II: Reziprozit\"atsgesetz, 
<Physica-Verlag, W\"urzburg-Vienna 1970 iv+204 pp.>
\newcount\hasse \global\hasse=\refno

\ref{\pc MILNE} (J):
Class field theory,
<http://www.jmilne.org/math.>
\newcount\milne \global\milne=\refno

\ref{\pc NEUKIRCH} (J):
Klassenk\"rpertheorie,
<Bibliographisches Institut, Mann\-heim, 1969. x+308.>
\newcount\neukirch \global\neukirch=\refno  
 
\ref{\pc SERRE} (J):
Corps locaux,
<Publications de l'Universit\'e de Nancago,
No.~VIII. Hermann, Paris, 1968. 245 pp.>
\newcount\serre \global\serre=\refno 

\ref{\pc TATE} (J):
Problem 9 : The general reciprocity law,
<in Mathematical developments arising from Hilbert problems
(Proc. Sympos. Pure Math., Northern Illinois Univ., De Kalb, Ill.,
1974), pp. 311--322. Proc. Sympos. Pure Math., Vol. XXVIII,
Amer. Math. Soc., Providence, R. I., 1976. > 
\newcount\tate \global\tate=\refno

\ref{\pc WEIL} (A):
Basic number theory.,
<Classics in Mathematics. Springer-Verlag, Berlin, 1995. xviii+315.>
\newcount\weil \global\weil=\refno

}

\centerline{\bf Classical Reciprocity Laws}

\vskip1cm
\centerline{Chandan Singh Dalawat}
\medskip

{\eightpoint Using the quadratic reciprocity law as the motivating
  example, we convey an understanding of classical reciprocity laws.}

\bigskip

Every since Gau\ss\ published his {\it Disquisitiones\/} in 1801,
reciprocity laws have been one of the main preoccupations of
arithmeticians.  My purpose here is not to go into the history of
these laws but to convey our present understanding.  The adjective
{\it classical\/} refers to the fact that we assume the relevant roots
of unity are present in the number field under discussion.

Let's begin with the quadratic reciprocity law.  For every prime
number $p$, we have the finite field $\F_p=\Z/p\Z$.  Its
multiplicative group $\F_p^\times$ is cyclic of order $p-1$.  If
$p\neq2$, then $p-1$ is even, so there is a unique surjective morphism
of groups $\lambda_p:\F_p^\times\to\Z^\times$, where $\Z^\times$ is
the multiplicative group consisting of $1$ and $-1$.

In general, $G$ being any group, a surjective morphism of groups
$G\to\Z^\times$ is called a {\it quadratic character\/} of $G$.  We've
seen that for every {\it odd\/} prime $p$, there is a {\it unique\/}
quadratic character of $\F_p^\times$.

For the even prime $2$, we need to consider the {\it three\/}
quadratic characters of the multiplicative group $(\Z/8\Z)^\times$.
The first one comes from the unique isomorphism of groups
$\lambda_4:(\Z/4\Z)^\times\to\Z^\times$~; indeed, the two groups in
question have order~$2$, so they are isomorphic and there is only one
isomorphism between them.  

To define the second one, view $\Z^\times$ as a subgroup of
$(\Z/8\Z)^\times$ (consisting of $\bar1$ and $-\bar1$), so that the
quotient group $(\Z/8\Z)^\times\!/\Z^\times$ has order~$2$.  There is
thus a unique isomorphism of groups
$\lambda_8:(\Z/8\Z)^\times\!/\Z^\times\to\Z^\times$, and it can be
viewed as a quadratic character of $(\Z/8\Z)^\times$.

The third quadratic character of $(\Z/8\Z)^\times$ is simply the
product $\lambda_4\lambda_8$, defined by
$\lambda_4\lambda_8(x)=\lambda_4(x)\lambda_8(x)$ for every
$x\in(\Z/8\Z)^\times$.  Out of these three quadratic characters, only
$\lambda_8$ is {\it even\/} in the sense that $\lambda_8(-1)=1$~; the
other two are {\it odd\/} because $\lambda_4(-1)=-1$ and
$\lambda_4\lambda_8(-1)=-1$. 

For every prime $p$, denote by $\Z_{(p)}$ the smallest subring of $\Q$
containg $l^{-1}$ for every prime $l\neq p$.  The morphism of rings
$\Z\to\F_p$ can be extended uniquely to a morphism of rings
$\Z_{(p)}\to\F_p$~; its kernel is $p\Z_{(p)}$.  We thus get a morphism
of groups $\Z_{(p)}^\times\to\F_p^\times$ which is easily seen to be
surjective.  For $p\neq2$, we may thus view $\lambda_p$ as a quadratic
character of $\Z_{(p)}^\times$.

Similarly, for $p=2$, the morphism of rings $\Z\to\Z/8\Z$ can be
extended uniquely to a morphism of rings $\Z_{(2)}\to\Z/8\Z$ (with
kernel $8\Z_{(2)}$).  We thus get a morphism of groups
$\Z_{(2)}^\times\to(\Z/8\Z)^\times$ which is easily seen to be
surjective.  We may thus view $\lambda_4$ and $\lambda_8$ as quadratic
characters of $\Z_{(2)}^\times$.

Till now, we've only defined some quadratic characters.  Here's our
first observation~: for every $a\in\Z_{(2)}^\times$,
$$
\lambda_4(a)=(-1)^{a-1\over2},\qquad
\lambda_8(a)=(-1)^{a^2-1\over8}.
$$ Without going into the proof, we clarify that these formulae have a
meaning.  Indeed, if $a\in\Z_{(2)}^\times$, then $a-1\in2\Z_{(2)}$, so
${a-1\over 2}\in\Z_{(2)}$.  Now, $\Z^\times$ is a (multiplicatively
written) vector space over $\F_2$ (of dimension~$1$), and hence a
module over $\Z_{(2)}$, so the expression $(-1)^{a-1\over2}$ has a
meaning.  Similarly, the expression $(-1)^{a^2-1\over8}$ has a meaning
for every $a\in\Z_{(2)}^\times$, for then $a^2-1\in8\Z_{(2)}$ and
${a^2-1\over8}\in\Z_{(2)}$.

The foregoing formulae can be said to compute $\lambda_4$ and
$\lambda_8$.  Can we compute $\lambda_p$ for odd primes $p$~?  In
other words, is there a formula for $\lambda_p(a)$, valid for every
$a\in\Z_{(p)}^\times$~?  Such a formula is precisely what the law of
quadratic reciprocity gives.

Recall that for $p\neq2$ the group $\Z_{(p)}^\times$ is collectively
generated by $-1$, $2$ and all odd primes $q\neq p$.  Since
$\lambda_p$ is a morphism of groups, it is sufficient to give a
formula for $\lambda_p(-1)$, $\lambda_p(2)$ and $\lambda_p(q)$.

The law of quadratic reciprocity asserts that for every prime $p\neq2$,
$$
\lambda_p(-1)=\lambda_4(p),\qquad
\lambda_p(2)=\lambda_8(p),\qquad
\lambda_p(q)=\lambda_q(\lambda_4(p)p),
$$ 
for every odd prime $q\neq p$.  It was first discovered by Euler and
Legendre in their old age, and proved by the young Gau\ss.  Since
then, a number of different proofs have been given (many of them by
Gau\ss\ himself), and it has been vastly generalised.

One the simplest proofs of the law
$\lambda_p(q)=\lambda_q(\lambda_4(p)p)$ is perhaps the one given by
Rousseau in 1991.  It consists in computing the product of all
elements in the group $(\F_p^\times\times\F_q^\times)/\Z^\times$ in
two different ways, by exploiting the isomorphism of groups
$(\Z/pq\Z)^\times\to\F_p^\times\times\F_q^\times$.

The law of quadratic reciprocity was generalised by Gau\ss, Jacobi,
and Eisenstein to cubic and quartic reciprocity laws.  For this
purpose, they had to enlarge the field $\Q$ to $\Q(j)$ and $\Q(i)$
respectively, where $j$ is a primitive third root of $1$ ($j^3=1$,
$j\neq1$) and $i$ is a primitive fourth root of $1$ ($i^4=1$,
$i^2\neq1$).  Dirichelet found the analogue of quadratic reciprocity
for the field $\Q(i)$.  Eisenstein and Kummer made deep contributions
to some cases of the $l$-tic reciprocity law in $\Q(\zeta)$, where
$\zeta$ is a primitive $l$-th root of unity ($\zeta^l=1$,
$\zeta\neq1$) and $l$ is an odd prime.

But let us jump directly to Hilbert, who reformulated the quadratic
reciprocity law as a {\it product formula\/} which made it possible to
guess what the generalisation to $m$-tic reciprocity should be, for
every $m>1$ (over a number field which contains a primitive $m$-th
root of unity).

The first notion we need is that of a {\it place\/} of $\Q$, which can
be either finite or archimedean.  A finite place of $\Q$ is just a
prime number, and there is just one archimedean place, denoted
$\infty$.  Shortly we shall define the completion $\Q_v$ of $\Q$ at a
place $v$.  It will turn out that $\Q_\infty$ is just the field $\R$
of real numbers.  For every finite place $p$ of $\Q$, Hensel defined a
new field called the field of $p$-adic numbers and denoted $\Q_p$.  It
is in terms of these fields that the mystery in the following
definitions will be clarified.

For any two numbers $a,b\in\Q^\times$ and every place $v$ of $\Q$,
define $(a,b)_v\in\Z^\times$ by the following explicit but opaque
rules.

If $v=\infty$, then $(a,b)_\infty=1$ precisely when $a>0$ or $b>0$, so
that $(a,b)_\infty=-1$ if $a<0$ and $b<0$.  Here the definition is not
so mysterious because $(a,b)_\infty=1$ precisely when the equation
$ax^2+by^2=1$ has a solution $x,y\in\R$.

Now let $v$ be a finite place of $\Q$, so that it is some prime number
$p$.  Note that every $x\in\Q^\times$ can be uniquely written as
$x=p^{v_p(x)}u_x$, with $v_p(x)\in\Z$ and $u_x\in\Z_{(p)}^\times$.
Let $a,b\in\Q^\times$, write
$$
a=p^{v_p(a)}u_a,\quad
b=p^{v_p(b)}u_b,\quad 
(v_p(a),v_p(b)\in\Z, u_a,u_b\in\Z_{(p)}^\times), 
$$ 
and  put
$$
t_{a,b}
=(-1)^{v_p(a)v_p(b)}a^{v_p(b)}b^{-v_p(a)}
=(-1)^{v_p(a)v_p(b)}u_a^{v_p(b)}u_b^{-v_p(a)},
$$ 
which is visibly in $\Z_{(p)}^\times$.   Define
$$
(a,b)_p=\lambda_p(t_{a,b})\quad(p\neq2),\qquad
(a,b)_2=(-1)^{{u_a-1\over2}{u_b-1\over2}}\lambda_8(t_{a,b}).
$$
As we've said, these definitions might seem unmotivated and contrived,
but their real meaning will come out once we've defined the fields
$\Q_p$.

We are ready to state Hilbert's product formula.  It says that for
$a,b\in\Q^\times$, we have $(a,b)_v=1$ for almost all places $v$ of
$\Q$, and 
$$
\prod_v(a,b)_v=1.
$$  
By unravelling the definitions, this product formula can be seen to be
equivalent to the quadratic reciprocity law.

For example, when $p$ and $q$ are distinct odd (positive) primes, the
definitions give $(p,q)_\infty=1$,
$$
(p,q)_2=(-1)^{{p-1\over2}{q-1\over2}},\quad
(p,q)_p=\lambda_p(q),\quad
(p,q)_q=\lambda_q(p)
$$
and $(p,q)_l=1$ for every odd prime $l$ different from $p$ and $q$.
So the product formula in this case becomes
$$
(-1)^{{p-1\over2}{q-1\over2}}\lambda_p(q)\lambda_q(p)=1,
$$
which is equivalent to $\lambda_p(q)=\lambda_q(\lambda_4(p)p)$ in view
of $\lambda_4(p)=(-1)^{p-1\over2}$ and
$\lambda_q(-1)=\lambda_4(q)=(-1)^{q-1\over2}$. 

The advantage of this reformulation of the quadratic reciprocity law
as a product formula is that it is so neat, compact, and memorable.
The disadvantage is that one doesn't quite understand where the symbol
$(a,b)_p\in\Z^\times$ (for primes $p$) comes from.  It can be properly
understood only in terms of Hensel's $p$-adic numbers, to which we now
turn.

It is best to first define the ring $\Z_p$ of $p$-adic integers.  It
is the ``inverse limit'' of the system of rings $\Z/p^n\Z$ and
homomorphisms 
$$
\varphi_n:\Z/p^{n+1}\Z\to\Z/p^n\Z
$$ 
(of reduction modulo $p^n$). Thus a $p$-adic integer $x\in\Z_p$ is a
system of elements $x=(x_n)_{n>0}$ such that $x_n\in\Z/p^n\Z$ and
$\varphi_n(x_{n+1})=x_n$.  Addition and multiplication are defined
componentwise.  It turns out that the ring $\Z_p$ is integral, and
$\Z$ can be identified with a subring of $\Z_p$.  The field $\Q_p$ is
defined as the field of fractions of $\Z_p$.

The ring $\Z_p$ carries a natural topology, the coarsest topology for
which all the projection morphisms $\Z_p\to\Z/p^n\Z$ are continuous.
It induces a topology on $\Q_p$ for which it is locally compact and
$\Q$ is a dense subset.

For $a,b\in\Q_p^\times$, one can define the symbol
$(a,b)_p\in\Z^\times$ to be $1$ if the equation $ax^2+by^2=1$ has a
solution $x,y\in\Q_p$, and $(a,b)_p=-1$ otherwise.  One can check that
for $a,b\in\Q^\times$, this new definition in terms of the solvability
of $ax^2+by^2=1$ in $\Q_p$ coincides with the previous explicit
definition in terms of the quadratic characters $\lambda_*$.

There is an even better way of understanding the symbols
$(a,b)_v\in\Z^\times$ (for $v$ a place of $\Q$ and
$a,b\in\Q_v^\times$).  Let $M_v=\Q_v(\sqrt{\Q_v^\times})$ be the
maximal abelian extension of $\Q_v$ of exponent~$2$.  It is a {\it
  minor miracle\/} that there is a unique isomorphism
$$
r_v:\Q_v^\times\!/\Q_v^{\times 2}\to\Gal(M_v|\Q_v)
$$
such that for every extension $L$ of $\Q_v$ in $M_v$, the kernel of
the composite map
$$
\rho_L:\Q_v^\times\to\Q_v^\times\!/\Q_v^{\times
  2}\to\Gal(M_v|\Q_v)\to\Gal(L|\Q_v)
$$
is equal to the image of the norm map $L^\times\to\Q_v^\times$.  For
$a,b\in\Q_v^\times$, take $L=\Q_v(\sqrt b)$~; then $\rho_L(a)(\sqrt
b)$ is either $\sqrt b$ or $-\sqrt b$, and $(a,b)_v\in\Z^\times$ is
precisely the sign such that 
$$
\rho_L(a)(\sqrt b)=(a,b)_v\sqrt b\qquad (L=\Q_v(\sqrt b)).
$$

More is true.  Let $K$ be any finite extension of $\Q_v$, let $m>0$ be
any integer, and let $M$ be the maximal abelian extension of $K$ of
exponent dividing~$m$.  If $K$ happens to contain a primitive $m$-th
root of unity, then $M=K(\root m\of{K^\times})$, by Kummer theory.  It
is a {\it minor miracle\/} that there is a unique isomorphism
$$
r_{m,K}:K^\times\!/K^{\times m}\to\Gal(M|K)
$$
such that for every extension $L$ of $K$ in $M$, the kernel of
the composite map
$$ 
\rho_L:K^\times\to K^\times\!/K^{\times m}\to\Gal(M|K)\to\Gal(L|K)
$$ 
is equal to the image of the norm map $N_{L|K}:L^\times\to K^\times$.
[At the finite places, there is an additional requirement which we've
  omitted because we haven't defined the relevant concepts.  In short,
  let $M_0\subset M$ denote the maximal unramified extension of $K$ in
  $M$~; the groups $K^\times/N_{M_0|K}(M_0^\times)$ and $\Gal(M_0|K)$
  are both cyclic (of order~$m$) with a canonical generator, and the
  requirement is that $\rho_{M_0}$ send the generator of
  $K^\times/N_{M_0|K}(M_0^\times)$ to the generator of $\Gal(M_0|K)$.
  This requirement is automatic if $m=2$.  ]

With this minor miracle in hand, one could go on to discuss the
general reciprocity law, but let's stick to the classical case where
the presence of a primitive $m$-th root of unity is required.

Suppose therefore that $K$ contains a primitive $m$-th root of unity.
For $a,b\in K^\times$, take $L=K(\root m\of b)$~; then
$\rho_L(a)(\root m\of b)$ and $\root m\of b$ differ by an $m$-th root
of unity in $\mu_m$, and one can define $(a,b)_{m,K}\in\mu_m$ by the
requirement that 
$$
\rho_L(a)(\root m\of b)=(a,b)_{m,K}\root m\of b\qquad (L=K(\root m\of b)).
$$ 
This is a generalisation of the previous case $m=2$, $K=\Q_v$, where
we denoted $(a,b)_{m,K}$ simply by $(a,b)_v$.  This is the local
ingredient we need in order to state the classical reciprocity laws.

Let us now turn to a finite extension $F$ of $\Q$ (also called a
number field) and explain what is meant by a {\it place\/} of $F$.  As
in the case of $\Q$ above, places come in two varieties~: finite and
archimedean.  To a finite place $v$ corresponds a prime number~$p$,
and $v$ is said to be a $p$-adic place. An archimedean place can be
real or imaginary.

A real place of $F$ is simply an embedding $F\to\R$.  An imaginary
place of $F$ is an embedding $\iota:F\to\C$ such that
$\iota(F)\not\subset\R$, except that two embeddings $\iota_1,\iota_2$
determine the same imaginary place if they differ by the conjugation
$z\mapsto\bar z$ ($i\mapsto-i$) in $\C$~: if
$\iota_1(a)=\overline{\iota_2(a)}$ for every $a\in F$.  We see that an
archimedean place of $F$ is really an equivalence class of embeddings
$F\to\C$, two embeddings being equivalent if they differ by an element
of $\Gal(\C|\R)$.  Every $F$ has at least one and at most finitely
many archimedean places.

Similarly, for every prime number $p$, a $p$-adic place of $F$ is an
equivalence class of embeddings $F\to\bar\Q_p$, where two embeddings
are equivalent if they differ by an element of $\Gal(\bar\Q_p|\Q_p)$.
Here, $\bar\Q_p$ is a fixed algebraic closure of $\Q_p$.  Every $F$
has at least one and at most finitely many $p$-adic places (for every
prime $p$).

Recall that the field $\Q_p$ ($p$ prime) carries a natural topology
which makes it a locally compact field.  As a result, every algebraic
closure of $\Q_p$ also carries a natural topology (but $\bar\Q_p$ is
not locally compact).  Also, the group $\Gal(\bar\Q_p|\Q_p)$ is far
more complicated than $\Gal(\C|\R)$.

The next notion we need is that of the completion $F_v$ of $F$ at a
place $v$. If $v$ is real, then $F_v=\R$.  If $v$ is imaginary, then
$F_v=\C$.  For a $p$-adic place $v$ of $F$, the completion $F_v$ is
defined to be the closure of $\iota(F)$ in $\bar\Q_p$, where
$\iota:F\to\bar\Q_p$ is an embedding representing the place $v$.  It
is a finite extension of $\Q_p$, uniquely determined by $F$ and $v$,
and $[F:\Q]=\sum_{v|p}[F_v:\Q_p]$, just as
$[F:\Q]=\sum_{v|\infty}[F_v:\R]$, where $v|p$ means that $v$ is a
$p$-adic place and $v|\infty$ means that $v$ is an archimedean place.

Till now the number field $F$ has been arbitrary.  We have defined the
notion of a place $v$ of $F$, and the completion $F_v$ of $F$ at $v$.
Now let $m>0$ be an integer, and suppose that $F$ contains a primitive
$m$-th root of unity.  [If $m>2$, then $F$ does not have any real
  places.]  Clearly, every completion $F_v$ also contains a primitive
$m$-th root of unity, making it possible to define
$(a,b)_{m,F_v}\in\mu_m$ for $a,b\in F_v^\times$, as we've seen.  [At
  an imaginary place $v$, we take $(a,b)_{m,\C}=1$ for all
  $a,b\in\C^\times$ because $\C$ is algebraically closed.  For this
  reason, imaginary places play no role in what follows.]

The main theorem, which encompasses all classical reciprocity laws,
states that if we start with $a,b\in F^\times$, then $(a,b)_{m,F_v}=1$
for almost all $v$, and the product formula
$$
\prod_v(a,b)_{m,F_v}=1
$$ 
holds (in the group $\mu_m$).  Quadratic reciprocity is the special
case $F=\Q$, $m=2$.  Cubic reciprocity (which we haven't recalled) is
the special case $F=\Q(j)$, $m=3$.  Quartic reciprocity is the special
case $F=\Q(i)$, $m=4$.  Dirichlet's analogue of quadratic reciprocity
is the special case $F=\Q(i)$, $m=2$.  Eisenstein, Kummer and Takagi's
work on $l$-tic reciprocity (for an odd prime $l$) is the special case
$F=\Q(\zeta)$ (where $\zeta^l=1$, $\zeta\neq1$) and $m=l$.  

When $F$ is an arbitrary number field and $m=2$, we get the quadratic
reciprocity law in $F$, due to Hilbert.  We have $(a,b)_{2,F_v}=1$ if
and only if the equation $ax^2+by^2=1$ has a solution $x,y\in F_v$~;
otherwise $(a,b)_{2,F_v}=-1$, just as in the special case $F=\Q$.

It is difficult to appreciate just how much information is packed into
this single neat product formula.  To unravel this information in the
case of some particular number field $F$ (containing a primitive
$m$-th root of unity), we need to determine the places of $F$, and
more importantly to give an explicit formula for $(a,b)_{m,F_v}$ at
every place $v$.  This quest has given rise to some of the deepest and
most sublime mathematics ever dreamt of by a human mind.

The only shortcoming of the above product formula is that it is
applicable only to those number fields which contain a primitive
$m$-th root of unity.  This restriction has been removed by Takagi,
Artin and Hasse, who came up with the general reciprocity law.  I hope
to discuss it on some future occasion and show how Chevalley's
invention of id{\`e}les provides a conceptual understanding of the
general law, just as Hensel's invention of $p$-adic numbers provides a
conceptual understanding of the classical laws.

Let us end with the provenance of the word {\it reciprocity\/}.  It
was first used by Legendre to reflect the fact that when $p$ and $q$
are distinct odd primes and one of them is $\equiv1\mod 4$, then
$\lambda_p(q)=\lambda_q(p)$, which is sometimes written more simply as
$(q/p)=(p/q)$.  In words~: the value of $\lambda_p$ at $q$ is the same
as the value of $\lambda_q$ at $p$, or $q$ is a square modulo~$p$ if
and only if $p$ is a square modulo~$q$.  The meaning of the word got
reinforced with similar formulae such as $(a/b)=(b/a)$ which express
special cases of other classical reciprocity laws.  Since then, the
notion of reciprocity has become a central tenet of Arithmetic,
largely thanks to Robert Langlands, as attested by Roger Godement :

{\it Legendre a devin{\'e} la formule et Gauss est devenu instatan{\'e}ment
c{\'e}l{\`e}bre en la prouvant.  En trouver des
g{\'e}n{\'e}ralisations, par exemple aux anneaux d'entiers
alg{\'e}briques, ou d'autres d{\'e}monstrations a constitu{\'e} un
sport national pour la dynastie allemande suscit{\'e} par Gauss
jusqu'{\`a} ce que le reste du monde, {\`a} commencer par le Japonais
Takagi en 1920 et {\`a} continuer par Chevalley une dizaine
d'ann{\'e}es plus tard, d{\'e}couvre le sujet et, apr{\`e}s 1945, le
fasse exploser.  Gouvern{\'e} par un Haut Commissariat qui surveille
rigoureusement l'alignement de ses Grandes Pyramides, c'est
aujourd'hui l'un des domaines les plus respect{\'e}s des
Math{\'e}matiques.} \citer\godement(p.~313).

\medbreak

The paper which has most influenced my point of view is
\citer\tate().  The reader who wants to see the first part of this Note 
worked out in every detail can consult my {\it Six lectures on
  quadratic reciprocity} \citer\dalawat().

The study of reciprocity laws led to class field theory. There is a
fairly large number of books on this subject, starting with
Hasse \citer\hasse(), Chevalley \citer\chevalley() and
Artin-Tate \citer\artintate().  A comprehensive account can be found
in \citer\casselsfrohlich().  Other sources include the books by
Weil \citer\weil(), Serre \citer\serre(), and Neukirch \citer\neukirch(),
and the online notes of Milne \citer\milne().

\medbreak

{\bf Acknowledgments.}  I heartily thank Sujatha for her invitation to
write this Note and Franz Lemmermeyer for a critical reading.

\bigbreak
\unvbox\bibbox

\bigskip\bigskip
{\obeylines\parskip=0pt\parindent=0pt
Chandan Singh Dalawat
Harish-Chandra Research Institute
Chhatnag Road, Jhunsi, {\pc ALLAHABAD} 211\thinspace019, India
\vskip5pt
\tt dalawat@gmail.com}

\bye